\documentclass[12pt]{amsart}
\usepackage{amssymb}
\usepackage{xcolor}
\usepackage{graphicx}
\usepackage{float}
\usepackage{dsfont}
\usepackage{url}
\usepackage{hyperref}
\usepackage{mathtools}
\mathtoolsset{showonlyrefs}
\usepackage{anysize}
\usepackage{geometry}
\marginsize{2cm}{2cm}{2cm}{2cm}
\usepackage[normalem]{ulem}
\usepackage[mathlines]{lineno}
\theoremstyle{plain}
\numberwithin{equation}{section}
\newtheorem{thm}{Theorem}[section]

\newcommand{\PP}[1]{\left(#1\right)}
\newcommand{\CC}[1]{\left[#1\right]}
\newcommand{\LL}[1]{\left\{#1\right\}}

\newcommand{\abs}[1]{\left|#1\right|}
\newcommand{\norm}[1]{\left|\left|#1\right|\right|}
\renewcommand{\Pr}[1]{\mathbb{P}\left(#1\right)}

\newcommand{\R}{\mathbb{R}}

\newcommand{\dd}{\textnormal{d}}

\title[Distribution of the angle between random segments]{Distribution of the angle between random segments}
\author{Paulo Manrique-Mir\'on}
\email{manriquemiron@gmail.com}
\begin{document}
\begin{abstract}
The study of \textit{random segments} is a classical problem in geometrical probability whose answer depends on the mechanism used to generate the segments. We consider four independent random points uniformly distributed in the unit disk and form the two labeled segments $S_{AB}$ and $S_{CD}$. The random variable of interest is the usual angle in $[0,\pi]$ between the vectors $B-A$ and $D-C$, conditional on the event that the two segments intersect. By introducing normal and tangential coordinates for each supporting line, we obtain an integral expression for the conditional density. The signed tangential coordinate of the intersection point is retained throughout the argument, which prevents the reflected-root overcounting that arises when only its squared norm is used. The same change of variables also recovers known distributions for the distance of a random chord from the center and for the length of a random segment. The resulting intersection probability is $\frac{1}{3}\left(1-\frac{35}{12\pi^2}\right)$.
\end{abstract}
\maketitle

\section{\textbf{Introduction}}\label{sec:Intro}

Geometrical probability deals with the study of classical geometry objects, points, segments, lines, planes, circles, spheres, etc., which are generated through some random mechanism \cite{kendall1963geometrical,mathai1999introduction}.

The development of this area dates back at least to 1733 when Georges-Louis Leclerc, Comte de Buffon, wrote ``Mémoire sur le jeu de franc-carreau'', where he proposed and solved (not always correctly) three problems formulated as mathematical games. These are the clean-tiles problem, the needle problem, and the mesh problem \cite{mathai1999introduction}. Among them, possibly the best known is the needle problem, which consists of randomly throwing a needle of length $l$ on a set of equidistant parallel lines, with separation $h$, on the plane. The question is to determine the probability that the needle cuts a line. In fact, the value of this probability is $\frac{2l}{\pi h}$. This game allows us to set up a method of simulation to determine the value of $\pi$.

The analysis of a random geometric object depends on the random mechanism used to generate it. For example, Kendall and Moran \cite{kendall1963geometrical} discuss Bertrand's problem: determine the probability that a \textit{random chord} of a circle is longer than the side of the equilateral triangle inscribed in it. Three different interpretations of a random chord are considered. In the first model, the chord is formed by joining two points generated independently and uniformly on the circumference. In the second, the chord is perpendicular to a fixed diameter and its intersection point with that diameter is uniformly distributed along the diameter. In the third, a point is chosen uniformly in the disk and the chord through that point is perpendicular to the corresponding radius. The resulting probabilities are $1/3$, $1/2$, and $1/4$, respectively \cite{garwood1966distance}.

Garwood and Holroyd \cite{garwood1966distance} interpret a \textit{random chord} as the segment passing through two independent and uniformly distributed points $P,Q$ in the interior of a circle with radius one. They computed the density function of the distance $L$ of the chord to the center of the circle,
\[
f_L(l) = \frac{16}{3\pi} (1-l^2)^{3/2} \mathds{1}_{\LL{l\in[0,1]}},
\]since this distance determines the length of the chord.

Previously, Garwood and Tanner \cite{garwood19582800} found the density of the distance $D$ between $P$ and $Q$,
\[
f_D(d) = \frac{2d}{\pi}\PP{2 \arccos\PP{\frac{d}{2}} - \sin\PP{2\arccos\PP{\frac{d}{2}} }} \mathds{1}_{\LL{d \in [0,2]}}.
\]

In both works the \textit{infinitesimal strategy} is used to determine the densities of considered lengths, which consists in the following idea: if $f(w)$ is the density of $W$ then, intuitively, $f(w)dw$ is the probability that $W\in[w,w+dw]$.

\begin{figure}[h]
\centering
\includegraphics[width=0.8\textwidth]{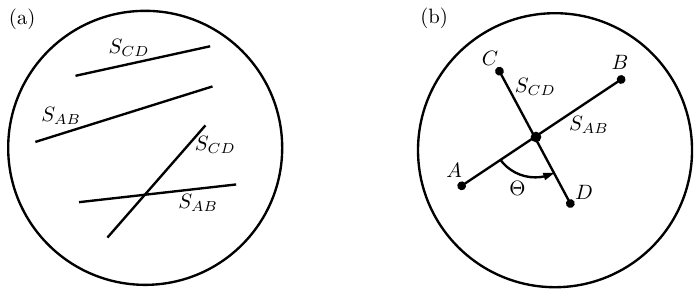}
\caption{}\label{fig:mie24may0940}
\end{figure}

In this manuscript, a \textit{random segment} is generated by joining two independent points uniformly distributed in the unit disk, following the mechanism used by Garwood and Holroyd. Two independent labeled segments are generated in this way, and we study the angle between their orientation vectors when the segments intersect; see Figure~\ref{fig:mie24may0940}.

Formally, let
\[
\mathds{D}=\{x\in\mathbb{R}^2:\|x\|\leq1\}.
\]
Let $X$ and $Y$ be independent random points, each uniformly distributed in $\mathds{D}$. The associated labeled random segment is
\[
S_{XY}:=\{(1-\alpha)X+\alpha Y:\alpha\in[0,1]\}.
\]
We consider four independent random points $A,B,C,D$, all uniformly distributed in $\mathds{D}$, and the segments $S_{AB}$ and $S_{CD}$. Since the endpoints are labeled, their orientation vectors are $B-A$ and $D-C$. Except on a null event, both vectors are nonzero. We define the usual angle between them by
\[
\Theta:=\arccos\left(\frac{(B-A)^T(D-C)}{\|B-A\|\,\|D-C\|}\right)\in[0,\pi].
\]
Our objective is to compute
\begin{equation}\label{eqn:mie24may0932}
\Pr{\Theta\leq\theta\mid S_{AB}\cap S_{CD}\neq\varnothing},
\end{equation}
for $\theta\in[0,\pi]$. 

To compute~\eqref{eqn:mie24may0932}, we introduce a change of variables that also recovers the results of Garwood--Holroyd and Garwood--Tanner.

The manuscript contains two further sections. Section~\ref{sec:MainResult} states the conditional density and records consequences of the change of variables. Section~\ref{sec:Proof} gives the proof.

\section{\textbf{Main Result}} \label{sec:MainResult}

The main result of this manuscript is presented below.
\begin{thm}\label{thm:mar29ago1553}
For $0<\theta<\pi$, define
\begin{align}
 g^*(\theta)
 :=\frac{1}{2\pi}\sum_{\varepsilon\in\{-1,1\}}
 \int_0^1\!\int_0^1
 &g^*_{1,\varepsilon}(\rho_{AB},\rho_{CD},\theta)\nonumber\\
 &\times
 \mathds{1}_{\left\{
 \sqrt{1-\rho_{AB}^2}\,|\sin\theta|
 \geq |\rho_{AB}\cos\theta+\varepsilon\rho_{CD}|
 \right\}}
 \,\dd\rho_{AB}\,\dd\rho_{CD},
 \label{eq:corrected-gstar}
\end{align}
where
\begin{align*}
 g^*_{1,\varepsilon}(\rho_{AB},\rho_{CD},\theta)
 &:=\left(\frac{4}{\pi}\right)^2
 \sqrt{1-\rho_{AB}^2}\sqrt{1-\rho_{CD}^2}\\
 &\quad\times
 \left[
 1-\frac{\rho_{AB}^2+\rho_{CD}^2
 +2\varepsilon\rho_{AB}\rho_{CD}\cos\theta}
 {\sin^2\theta}
 \right]^2
 \mathds{1}_{\{\rho_{AB},\rho_{CD}\in[0,1]\}}.
\end{align*}
Set $g^*(0)=g^*(\pi)=0$ and
\[
 c:=\int_0^\pi g^*(s)\,\dd s.
\]
Then
\[
 c=\Pr{S_{AB}\cap S_{CD}\neq\varnothing}
 =\frac13\left(1-\frac{35}{12\pi^2}\right)
 =0.234826627014\ldots,
\]
and the conditional density of $\Theta$ given intersection is
\[
 g(\theta)=\frac{g^*(\theta)}{c}\,
 \mathds{1}_{\{\theta\in(0,\pi)\}}.
\]
Equivalently,
\[
\Pr{\Theta\leq\theta\mid S_{AB}\cap S_{CD}\neq\varnothing}
 =\int_0^\theta g(s)\,\dd s,
 \qquad 0\leq\theta\leq\pi.
\]
\end{thm}

\begin{figure}[h]
\centering
\includegraphics[width=0.7\textwidth]{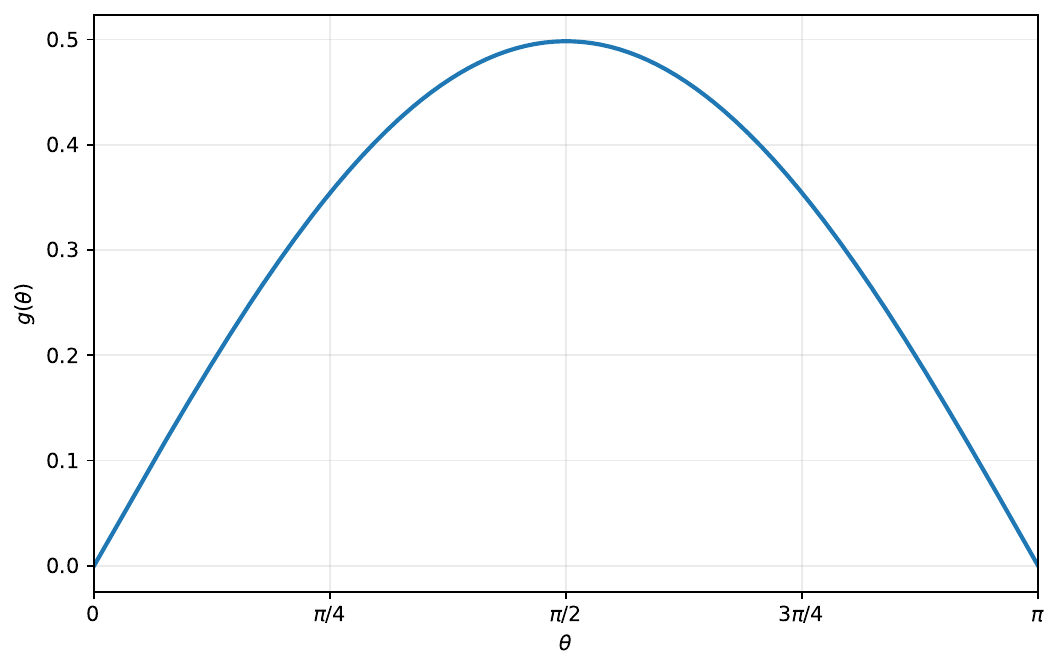}
\caption{Corrected conditional density of the usual angle between the labeled segment vectors.}\label{fig:vie25agosto1218}
\end{figure}

Figure~\ref{fig:vie25agosto1218} shows the conditional density $g(\theta)$. Its symmetry about $\pi/2$ follows from exchanging the labels $C$ and $D$.

The proof of Theorem ~\ref{thm:mar29ago1553} is based on the following change of variable.

Observe that $X=\sqrt{R_X}(\cos\Gamma_X,\sin\Gamma_X)^T$, where $v^T$ denotes transpose, is uniformly distributed in $\mathds{D}$ when $R_X$ and $\Gamma_X$ are independent, $R_X\sim\operatorname{Unif}[0,1]$, and $\Gamma_X\sim\operatorname{Unif}[0,2\pi]$. Here $R_X$ is the squared radial coordinate. Thus, $S_{AB}$ can be written as
\[
S_{AB} = \LL{w\in\R^2 : w=(1-\alpha) \sqrt{R_A}\begin{pmatrix}
\cos(\Gamma_A)\\
\sin(\Gamma_A)
\end{pmatrix}
+ \alpha \sqrt{R_B}\begin{pmatrix}
\cos(\Gamma_B)\\
\sin(\Gamma_B)
\end{pmatrix}, \; \alpha\in[0,1]},
\]

\noindent where $R_A,R_B,\Gamma_A,\Gamma_B$ are independent, $R_A,R_B\sim\operatorname{Unif}[0,1]$, and $\Gamma_A,\Gamma_B\sim\operatorname{Unif}[0,2\pi]$.

\begin{figure}[h]
\centering
\includegraphics[width=0.5\textwidth]{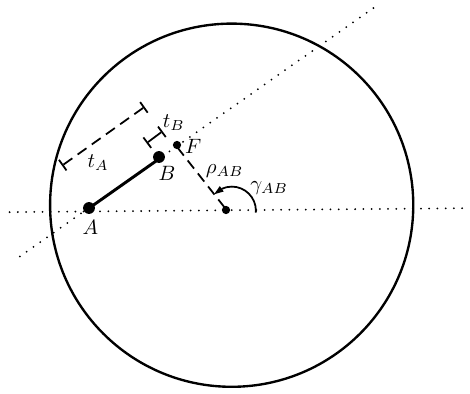}
\caption{}\label{fig:mie24may1115}
\end{figure}

Consider the perpendicular from the origin $O$ to the supporting line of $S_{AB}$, and denote its foot by $F_{AB}$. Let $\Gamma_{AB}$ be the angle made by $OF_{AB}$ with the $x$-axis and let $R_{AB}$ be the perpendicular distance from the origin to the supporting line. Define
\[
 n(\gamma):=(\cos\gamma,\sin\gamma)^T,
 \qquad
 u(\gamma):=(-\sin\gamma,\cos\gamma)^T.
\]
Then the endpoints are represented as
\[
 A=R_{AB}n(\Gamma_{AB})+T_Au(\Gamma_{AB}),
 \qquad
 B=R_{AB}n(\Gamma_{AB})+T_Bu(\Gamma_{AB}).
\]
Thus $|T_A|$ and $|T_B|$ are the distances of $A$ and $B$ from $F_{AB}$, respectively. We denote realizations of $R_{AB},\Gamma_{AB},T_A,T_B$ by $\rho_{AB},\gamma_{AB},t_A,t_B$. See Figure~\ref{fig:mie24may1115}.

Note that
\begin{align}\label{eqn:mie24may1126}
\sqrt{\rho_j} \cos\gamma_j & = \rho_{AB} \cos\gamma_{AB} - t_{j}\sin\gamma_{AB},\\
\sqrt{\rho_j} \sin\gamma_j & = \rho_{AB} \sin\gamma_{AB} + t_{j}\cos\gamma_{AB},\nonumber 
\end{align} 
for $j\in\LL{A,B}$.

Thus, the joint density of $(R_A,\Gamma_A,R_B,\Gamma_B)$ is
\[
f(\rho_A,\gamma_A,\rho_B,\gamma_B)= \frac{1}{(2\pi)^2} \mathds{1}_{\LL{\rho_A\in[0,1]}} \mathds{1}_{\LL{\gamma_A\in [0,2\pi]}} \mathds{1}_{\LL{\rho_B\in [0,1]}} \mathds{1}_{\LL{\gamma_B\in [0,2\pi]}},
 \]
and can be expressed in terms of $\rho_{AB}, \gamma_{AB}, t_A, t_B$ as
\begin{align}\label{eqn:mie24may1226}
& \PP{\frac{2}{2\pi}}^2\abs{t_A-t_B} \mathds{1}_{\LL{\rho_{AB}\in[0,1]}} \mathds{1}_{\LL{\gamma_{AB}\in[0,2\pi]}} \mathds{1}_{\LL{t_A\in\CC{-\sqrt{1-\rho^2_{AB}},\sqrt{1-\rho^2_{AB}}}}} \mathds{1}_{\LL{t_B\in\CC{-\sqrt{1-\rho^2_{AB}},\sqrt{1-\rho^2_{AB}}}}}. 
\end{align}

This change of variable allows us to obtain more information of the random segment that does not seem clear from its original definition. For example, the results of Garwood and Holroyd and Garwood and Tanner can be deduced directly from this. The marginal density at $\rho_{AB}$ retrieves the result from Garwood and Holroyd,
\begin{align*}
f(\rho_{AB}) & = \mathds{1}_{\LL{\rho_{AB}\in[0,1]}} \frac{2}{\pi} \int \int \abs{t_A-t_B} \mathds{1}_{\LL{t_A\in\CC{-\sqrt{1-\rho^2_{AB}},\sqrt{1-\rho^2_{AB}}}}} \mathds{1}_{\LL{t_B\in\CC{-\sqrt{1-\rho^2_{AB}},\sqrt{1-\rho^2_{AB}}}}} \dd t_A \dd t_B \\
& = \frac{16}{3\pi} (1-\rho_{AB}^2)^{3/2}  \mathds{1}_{\LL{\rho_{AB}\in[0,1]}}.
\end{align*}
Meanwhile, the marginal density at $(t_A,t_B)$
\begin{align*}
f(t_A,t_B) & = \frac{2}{\pi} \abs{t_A-t_B} \min\LL{\sqrt{1-t_A^2},\sqrt{1-t_B^2}} \mathds{1}_{\LL{t_A\in[-1,1]}} \mathds{1}_{\LL{t_B\in[-1,1]}},
\end{align*}
allows to retrieve the result of Garwood and Tanner,
\begin{align*}
\Pr{\abs{S_{AB}} \leq d} & = \Pr{\abs{T_A-T_B} \leq d} \\
& = \int_{\LL{\abs{t_A-t_B} \leq d}} f(t_A,t_B) \dd(t_A,t_B) \\
& = \frac{8}{\pi}\int_0^{d/2}\int_{-t_B}^{t_B}
(t_B-t_A)\sqrt{1-t_B^2}\,\dd t_A\,\dd t_B \\
&\quad +\frac{8}{\pi}\int_{d/2}^1\int_{t_B-d}^{t_B}
(t_B-t_A)\sqrt{1-t_B^2}\,\dd t_A\,\dd t_B \\
& = \int_0^d \frac{s}{\pi} \frac{-4s+s^3+8\sqrt{4-s^2}\,\textnormal{arccot}\PP{\frac{2+s}{\sqrt{4-s^2}}}}{\sqrt{4-s^2}} \mathds{1}_{\LL{s\in[0,2]}}  \dd s,
\end{align*} with a little extra algebraic work.

In the main result, the intersection condition creates dependence among $\Gamma_{AB},\Gamma_{CD},R_{AB}$, and $R_{CD}$, which prevents an elementary closed form for $g(\theta)$. Nevertheless, the integral representation gives a direct numerical scheme; see Figure~\ref{fig:vie25agosto1218}. Moreover,
\[
\Pr{S_{AB}\cap S_{CD}\neq\varnothing}
=\int_0^\pi g^*(\theta)\,\dd\theta
=\frac13\left(1-\frac{35}{12\pi^2}\right)
\approx0.2348266270.
\]
Thus, two independently generated segments intersect with probability about $23.48\%$. 

An appropriate change of variables makes the geometry and the source of dependence transparent. The remaining integral complexity is intrinsic to conditioning on the intersection of the two random segments.

\section{\textbf{Proof}} \label{sec:Proof}

In this section the proof of Theorem ~\ref{thm:mar29ago1553} is presented.

From~\eqref{eqn:mie24may1126}, for $j\in\{A,B\}$,
\begin{align}\label{eqn:mie24may1128}
\rho_j &= \rho_{AB}^2+t_j^2,\\
\cos\gamma_j
&=\frac{\rho_{AB}\cos\gamma_{AB}-t_j\sin\gamma_{AB}}
{\sqrt{\rho_{AB}^2+t_j^2}},\nonumber\\
\sin\gamma_j
&=\frac{\rho_{AB}\sin\gamma_{AB}+t_j\cos\gamma_{AB}}
{\sqrt{\rho_{AB}^2+t_j^2}},\nonumber\\
\gamma_j
&=\operatorname{atan2}\!\left(
\rho_{AB}\sin\gamma_{AB}+t_j\cos\gamma_{AB},
\rho_{AB}\cos\gamma_{AB}-t_j\sin\gamma_{AB}
\right).\nonumber
\end{align}
The two-argument function $\operatorname{atan2}$ is required to retain the correct quadrant; a one-argument arctangent is only defined modulo $\pi$. The two-argument arctangent function, denoted by
\[
\operatorname{atan2}(y,x),
\]
assigns to each nonzero vector
\((x,y)\in\mathbb{R}^{2}\setminus\{(0,0)\}\)
its polar angle, taking into account the signs of both \(x\) and \(y\).
More precisely, if
\[
r=\sqrt{x^{2}+y^{2}},
\]
then \(\theta=\operatorname{atan2}(y,x)\) is the unique angle
\(\theta\in(-\pi,\pi]\) satisfying
\[
\cos\theta=\frac{x}{r},
\qquad
\sin\theta=\frac{y}{r}.
\]
Unlike the function \(\arctan(y/x)\), the function
\(\operatorname{atan2}(y,x)\) correctly determines the quadrant
containing the point \((x,y)\) and remains well defined when \(x=0\),
provided that \(y\neq0\). The expression
\(\operatorname{atan2}(0,0)\) is undefined.

The joint density of $(R_A,\Gamma_A,R_B,\Gamma_B)$
\[
f(\rho_A,\gamma_A,\rho_B,\gamma_B)= \frac{1}{(2\pi)^2} \mathds{1}_{\LL{\rho_A\in[0,1]}} \mathds{1}_{\LL{\gamma_A\in [0,2\pi]}} \mathds{1}_{\LL{\rho_B\in [0,1]}} \mathds{1}_{\LL{\gamma_B\in [0,2\pi]}}
 \]
is written in terms of  $\rho_{AB}, \gamma_{AB}, t_{A}, t_{B}$, i.e.,
\begin{align*}
\hspace{-0.65cm} & f(\rho_{AB}, \gamma_{AB}, t_{A}, t_{B}) \\
\hspace{-0.65cm} & = \frac{1}{(2\pi)^2} \abs{J} \mathds{1}_{\LL{\rho_{AB}\in[0,1]}} \mathds{1}_{\LL{\gamma_{AB}\in[0,2\pi]}} \mathds{1}_{\LL{t_A\in\CC{-\sqrt{1-\rho^2_{AB}},\sqrt{1-\rho^2_{AB}}}}} \mathds{1}_{\LL{t_B\in\CC{-\sqrt{1-\rho^2_{AB}},\sqrt{1-\rho^2_{AB}}}}},
\end{align*}
where $\abs{J}$ is the absolute value of the determinant of the Jacobian matrix $J$, namely
\[
\renewcommand*{\arraystretch}{1.5}
J = \begin{pmatrix}
\frac{\partial \rho_A}{\partial t_A}  & \frac{\partial \rho_A}{\partial t_B} & \frac{\partial \rho_A}{\partial \rho_{AB}} & \frac{\partial \rho_A}{\partial \gamma_{AB}}\\

\frac{\partial \gamma_A}{\partial t_A}  & \frac{\partial \gamma_A}{\partial t_B} & \frac{\partial \gamma_A}{\partial \rho_{AB}} & \frac{\partial \gamma_A}{\partial \gamma_{AB}}\\

\frac{\partial \rho_B}{\partial t_A}  & \frac{\partial \rho_B}{\partial t_B} & \frac{\partial \rho_B}{\partial \rho_{AB}} & \frac{\partial \rho_B}{\partial \gamma_{AB}}\\

\frac{\partial \gamma_B}{\partial t_A}  & \frac{\partial \gamma_B}{\partial t_B} & \frac{\partial \gamma_B}{\partial \rho_{AB}} & \frac{\partial \gamma_B}{\partial \gamma_{AB}}

\end{pmatrix}
=
 \begin{pmatrix}
2t_A & 0 & 2\rho_{AB} & 0 \\
\frac{\rho_{AB}}{\rho_{AB}^2+t_A^2} & 0 & -\frac{t_A}{\rho_{AB}^2+t_A^2} & 1 \\
0 & 2t_B & 2\rho_{AB} & 0 \\
0 & \frac{\rho_{AB}}{\rho_{AB}^2+t_B^2} & -\frac{t_B}{\rho_{AB}^2+t_B^2} & 1
\end{pmatrix}.
\]
Then $\abs{J} = 4\abs{t_A-t_B}$. Thus, the joint density $f(\rho_{AB}, \gamma_{AB}, t_{A}, t_{B})$ is
\begin{align}
& \PP{\frac{2}{2\pi}}^2\abs{t_A-t_B} \mathds{1}_{\LL{\rho_{AB}\in[0,1]}} \mathds{1}_{\LL{\gamma_{AB}\in[0,2\pi]}} \mathds{1}_{\LL{t_A\in\CC{-\sqrt{1-\rho^2_{AB}},\sqrt{1-\rho^2_{AB}}}}} \mathds{1}_{\LL{t_B\in\CC{-\sqrt{1-\rho^2_{AB}},\sqrt{1-\rho^2_{AB}}}}}. 
\end{align}

\begin{figure}[h]
\centering
\includegraphics{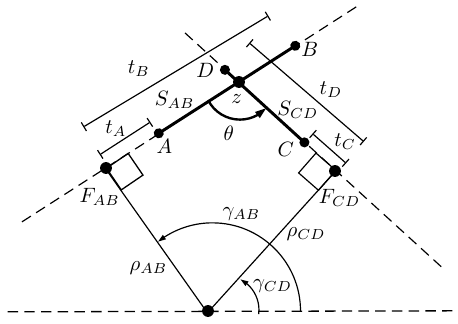}
\caption{}\label{fig:vie26may1213}
\end{figure}

From this point on, consider both segments in normal--tangential coordinates. For $S_{CD}$, use the analogous variables $\rho_{CD},\gamma_{CD},t_C,t_D$. The supporting lines are
\[
\ell_{AB}=\{x:n(\gamma_{AB})^Tx=\rho_{AB}\},
\qquad
\ell_{CD}=\{x:n(\gamma_{CD})^Tx=\rho_{CD}\}.
\]
Put $\delta:=\gamma_{AB}-\gamma_{CD}$. Since the angular variables have continuous distributions,
\[
\Pr{\sin\delta=0}=0.
\]
Consequently, the supporting lines have a unique intersection point almost surely. It is the solution of the nonsingular normal system
\begin{align}\label{eqn:vie26may1244}
 n(\gamma_{AB})^Tz&=\rho_{AB},\\
 n(\gamma_{CD})^Tz&=\rho_{CD}.\nonumber
\end{align}
Solving this system gives
\begin{equation}\label{eqn:vie26may1246}
 z=\frac1{\sin\delta}
 \begin{pmatrix}
 \rho_{CD}\sin\gamma_{AB}-\rho_{AB}\sin\gamma_{CD}\\
 \rho_{AB}\cos\gamma_{CD}-\rho_{CD}\cos\gamma_{AB}
 \end{pmatrix},
\end{equation}
and
\begin{equation}\label{eqn:vie26may1324}
 \|z\|^2
 =\frac{\rho_{AB}^2+\rho_{CD}^2
 -2\rho_{AB}\rho_{CD}\cos\delta}{\sin^2\delta}.
\end{equation}
The phrase ``always a unique solution'' must therefore be read as ``a unique solution almost surely''; parallel supporting lines form a null event. (See Figure~\ref{fig:vie26may1213}.)
 
The signed tangential coordinates of $z$ on the two lines are
\begin{align}
 q_{AB}&:=u(\gamma_{AB})^Tz
 =\frac{\rho_{AB}\cos\delta-\rho_{CD}}{\sin\delta},
 \label{eq:qAB}\\
 q_{CD}&:=u(\gamma_{CD})^Tz
 =\frac{\rho_{AB}-\rho_{CD}\cos\delta}{\sin\delta}.
 \label{eq:qCD}
\end{align}
These signs cannot be discarded. Indeed,
\[
 z\in S_{AB}
 \quad\Longleftrightarrow\quad
 (t_A-q_{AB})(t_B-q_{AB})\leq0,
\]
and similarly
\[
 z\in S_{CD}
 \quad\Longleftrightarrow\quad
 (t_C-q_{CD})(t_D-q_{CD})\leq0.
\]
Except on the null events $t_A=t_B$ or $t_C=t_D$, the corresponding affine parameters are unique:
\begin{align}
 \alpha_{AB}&=\frac{q_{AB}-t_A}{t_B-t_A},
 \label{eqn:dom28may0731}\\
 \alpha_{CD}&=\frac{q_{CD}-t_C}{t_D-t_C}.
 \label{eqn:dom28may0732}
\end{align}
Thus
\[
 S_{AB}\cap S_{CD}\neq\varnothing
 \quad\Longleftrightarrow\quad
 \mathcal I_{AB}\cap\mathcal I_{CD},
\]
where
\[
 \mathcal I_{AB}:=\{(t_A-q_{AB})(t_B-q_{AB})\leq0\},
 \qquad
 \mathcal I_{CD}:=\{(t_C-q_{CD})(t_D-q_{CD})\leq0\}.
\]
This exact event replaces the union over the two reflected roots obtained from the scalar norm equation \eqref{eqn:vie26may1324}. Squaring the tangential coordinate produces the artificial alternatives $\pm\sqrt{\|z\|^2-\rho^2}$; only the signed value $q$ represents the actual point $z$.

We now integrate the endpoint coordinates without overcounting. Fix $\rho\in[0,1]$ and let
\[
 a:=\sqrt{1-\rho^2},\qquad |q|\leq a.
\]
For the orientation $t_A<t_B$, the contribution of one segment is
\begin{align*}
 I_+(\rho,q)
 &=\frac1{\pi^2}\int_{-a}^{q}\int_q^a(y-x)\,\dd y\,\dd x\\
 &=\frac{a}{\pi^2}(a^2-q^2).
\end{align*}
By exchanging the two endpoints,
\[
 I_-(\rho,q)=I_+(\rho,q).
\]
Therefore the total contribution of one segment is
\begin{equation}\label{eq:one-segment-correct}
 I(\rho,q)=\frac{2a}{\pi^2}(a^2-q^2)
 =\frac{2}{\pi^2}\sqrt{1-\rho^2}\,(1-\|z\|^2),
\end{equation}
where $a^2-q^2=1-\|z\|^2$. This integral counts the actual signed coordinate $q$ once. Integrating both $q$ and $-q$ would double the contribution for each segment and would create a factor $2\times2=4$.

Define the orientation signs
\[
 \sigma_{AB}:=\operatorname{sgn}(t_B-t_A),
 \qquad
 \sigma_{CD}:=\operatorname{sgn}(t_D-t_C),
 \qquad
 \eta:=\sigma_{AB}\sigma_{CD}\in\{-1,1\}.
\]
For each fixed value of $\eta$, two of the four orientation combinations are possible. Hence, after integrating $t_A,t_B,t_C,t_D$ under the exact intersection event and a fixed product $\eta$, the contribution is
\begin{align}\label{eq:f3-correct-sign}
 f^{**}_{3,\eta}
 &=\frac{2}{\pi^4}
 \sqrt{1-\rho_{AB}^2}\sqrt{1-\rho_{CD}^2}
 (1-\|z\|^2)^2,
\end{align}
with the natural indicators $\rho_{AB},\rho_{CD}\in[0,1]$, $\gamma_{AB},\gamma_{CD}\in[0,2\pi]$, and $\|z\|\leq1$. Summing over $\eta=\pm1$ gives $4/\pi^4$ times the same geometric factor.

The angle between the segment vectors satisfies
\begin{equation}\label{eq:angle-sign}
 \cos\Theta
 =\eta\cos\delta,
\end{equation}
because
\[
 B-A=(t_B-t_A)u(\gamma_{AB}),
 \qquad
 D-C=(t_D-t_C)u(\gamma_{CD}).
\]
Let
\[
 \Phi:=\arccos(\cos\delta)\in[0,\pi].
\]
Then $\Theta=\Phi$ when $\eta=1$, whereas $\Theta=\pi-\Phi$ when $\eta=-1$.

The linear difference $\Gamma:=\Gamma_{AB}-\Gamma_{CD}$ has density
\[
 h(\gamma)=\frac{2\pi-|\gamma|}{4\pi^2}
 \mathds{1}_{\{\gamma\in[-2\pi,2\pi]\}}.
\]
For $0<\theta<\pi$, the four preimages of the line angle $\Phi=\theta$ are
$\pm\theta$ and $\pm(2\pi-\theta)$, and
\begin{equation}\label{eq:angular-folding}
 h(\theta)+h(-\theta)+h(2\pi-\theta)+h(-(2\pi-\theta))
 =\frac1\pi.
\end{equation}
The same total weight is obtained for the preimages of $\Phi=\pi-\theta$.

For a fixed $\theta$, write $\varepsilon=-\eta$. In the branch contributing to $\Theta=\theta$, formula~\eqref{eqn:vie26may1324} becomes
\begin{equation}\label{eq:z-epsilon}
 \|z\|^2
 =\frac{\rho_{AB}^2+\rho_{CD}^2
 +2\varepsilon\rho_{AB}\rho_{CD}\cos\theta}
 {\sin^2\theta}.
\end{equation}
Moreover, $\|z\|\leq1$ is equivalent to
\begin{equation}\label{eq:support-epsilon}
 \sqrt{1-\rho_{AB}^2}\,|\sin\theta|
 \geq|\rho_{AB}\cos\theta+\varepsilon\rho_{CD}|.
\end{equation}
Combining~\eqref{eq:f3-correct-sign}, the four angular preimages in~\eqref{eq:angular-folding}, and the factor $4\pi^2$ that converts the density of the two original angles to the density $h$ of their difference, the contribution of each branch is
\[
 \frac{1}{2\pi}\left(\frac4\pi\right)^2
 \sqrt{1-\rho_{AB}^2}\sqrt{1-\rho_{CD}^2}
 \left(1-\|z\|^2\right)^2.
\]
After substituting~\eqref{eq:z-epsilon}, imposing~\eqref{eq:support-epsilon}, integrating over $\rho_{AB},\rho_{CD}$, and summing over $\varepsilon\in\{-1,1\}$, we obtain exactly~\eqref{eq:corrected-gstar}. Therefore
\begin{equation}\label{eqn:vie30jun1100}
 \Pr{\Theta\leq\theta,\ S_{AB}\cap S_{CD}\neq\varnothing}
 =\int_0^\theta g^*(s)\,\dd s.
\end{equation}
Taking $\theta=\pi$ gives
\[
 c=\int_0^\pi g^*(s)\,\dd s
 =\Pr{S_{AB}\cap S_{CD}\neq\varnothing}.
\]
Finally, four points in general position admit three pairings. Exactly one pairing crosses when the four points are in convex position, and no pairing crosses when one point lies inside the triangle formed by the other three. By symmetry of the labels,
\[
 c=\frac13\Pr{A,B,C,D\text{ are in convex position}}.
\]

Thus,
\[
\Pr{A,B,C,D\text{ are in convex position}} = 1 - 4 \Pr{D\in \triangle(ABC)}.
\]
According our of change of variable, the area of $\triangle(ABC)$ is 
\[
\frac{1}{2}\abs{t_A-t_B}\abs{\rho_{AB} - u(\gamma_{AB})^T C}.
\]
Hence
\begin{align*}
\Pr{D\in \triangle(ABC)} & = \int \int \frac{\mbox{Area}(\triangle(ABC))}{\pi} \dd(A,B) \dd C \\
& =  \int \int \frac{1}{2\pi}\abs{t_A-t_B}\abs{\rho_{AB} - u(\gamma_{AB})^T C} \dd(T_A,T_B,R_{AB},\Gamma_{AB}) \dd C \\
& = \int \frac{\norm{C}^6-6\norm{C}^4+24\norm{C}^2+16}{36\pi^2} \dd C \\
& = \frac{35}{48\pi^2}.
\end{align*}
Then
\[
c= \frac{1}{3}\PP{ 1 - \frac{35}{12\pi^2}}.
\]

For four independent points uniformly distributed in the disk, the Sylvester probability is
\[
 \Pr{A,B,C,D\text{ are in convex position}}
 =1-\frac{35}{12\pi^2}
\]
(see, for example, \cite{kabluchko2026refinement}). This proves the stated value of $c$. Dividing~\eqref{eqn:vie30jun1100} by $c$ completes the proof of Theorem~\ref{thm:mar29ago1553}.



\end{document}